\documentclass[11pt,a4paper, final, twoside]{article}
\author{ Traian Preda}
\title{ A note of Zuk's citerion}
\usepackage{amsmath}
\usepackage{fancyhdr}
\usepackage{amsthm}
\usepackage{amsfonts}
\usepackage{amssymb}
\usepackage{amscd}
\usepackage{amsthm}
\usepackage{graphicx}
\usepackage{afterpage}
\usepackage[colorlinks=true, urlcolor=blue,  linkcolor=blue, citecolor=blue]{hyperref}

\newtheorem{theorem}{Theorem}

\newtheorem{definition}[theorem]{Definition}
\newtheorem{example}[theorem]{Example}

\numberwithin{equation}{section}

\begin{document}
\hyphenpenalty=10000
 
\begin{center}
{\Large \textbf{A note of Zuk's criterion }}\\[5mm]
{\large {Traian Preda  }\\[10mm]
}
\end{center}

{\footnotesize \textbf{Abstract}. Zuk's criterion give us a condition for a finitely 
generated group to have Property(T): the 
smallest
 non - zero eigenvalue of Laplace operator
 \( \Delta_\mu \)
 corresponding to the simple random walk on $\mathcal{G}$(S) satisfies   $\lambda_1(\mathcal{G}))>\frac{1}{2}$.
We present here two examples that prove that this condition cannot be improved.
}
 
\afterpage{
\fancyhead{} \fancyfoot{} 
\fancyhead[LE, RO]{\bf\thepage}
\fancyhead[LO]{\small A note of Zuk's criterion}
\fancyhead[RE]{\small Traian Preda  }
}

\begin{definition}(see \cite{bhv} and \cite{gro} )

i) A random walk or Markov kernel on a non-empty set X is a kernel with non-negative
values \(\mu : X\times X \to\mathbb{R}_+ \) such that:
\[
\sum_{y\in X}\mu (x,y) = 1, \forall x\in X.
\]

ii) A stationary measure for a random walk $\mu$ is a function \(\nu : X\to\mathbb{R}^*_+ \) such that :
\[
\nu (x) \mu (x, y) = \nu (y) \mu (y,x), \forall x,y \in X.
\]
\end{definition}

\begin{example}
Let \(\mathcal{G}\) =(X,E) be a locally finite graph.
For x,y $\in$ X, set
\begin{equation}
 \mu (x,y) = \left\{
\begin  {matrix}
{\displaystyle\frac{1}{ deg(x)} }& \quad if (x,y) \in E \\[0.7em]
0 & \quad otherwise 
\end {matrix}\right.
\end{equation}
and deg(x) =card \(\{y\in X| (x,y) \in E\} \) is the degree of a vertex x \(\in\) X.\\
$\mu$ is called simple random walk on X and $\nu$ is a stationary measure for $\mu$.
\end{example}

Consider the Hilbert space:
\[
\Omega^0_{\mathbb{C}}(X) = \{ f : X\to\mathbb{C} | \displaystyle\sum_{x\in X} |f(x)|^2\nu (x) < \infty \}
\]

The Laplace operator $\Delta_\mu$ on $\Omega^0_{\mathbb{C}}(X) $ is definited by
\(
(\Delta_\mu f)(x) = f(x) - \displaystyle\sum_{x\sim y}f(y)\mu(x,y)
\).

Let $\Gamma$ be a group generated by a finite set S. We assume that e $\notin$ S and  S = S$^{-1}$  (S~is~simetric).

The graph $\mathcal{G}(S)$ associated to $S$ has vertex set  $S$ and the set of edges is
the set of pairs $(s,t) \in S\times S $ such that $s^{-1}t \in S$.

\begin{theorem}(Zuk's criterion)( see \cite{zuk})

Let $\Gamma$ be a group generated by a finite set $S$ with $e \notin S$. Let $\mathcal{G}(S)$ be the graph associated to $S$. Assume that $\mathcal{G}(S)$ is connected and that the smallest
non-zero eigenvalue of the Laplace operator $\Delta_\mu$ corresponding to the simple random walk on $\mathcal{G}(S)$ satisfies $ \lambda_1(\mathcal{G}(S))>\displaystyle\frac{1}{2}$.

Then $\Gamma$ has Property (T). 
\end{theorem}
We prove that the condition $ \lambda_1(\mathcal{G}(S))>\displaystyle\frac{1}{2}$ cannot be improved, using two examples.
\begin{example}
Consider $S$ = \{ $1, -1, 2, -2$\} a generating set of the group $\mathbb{Z}$ and let
$\mathcal{G}(S)$ be the finite graph associated to $S$. Then the graph $\mathcal{G}(S)$ is
the graph:
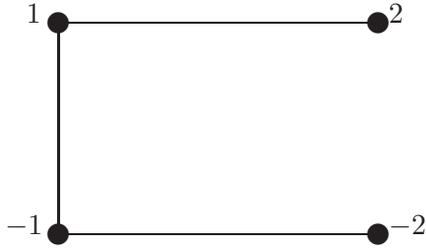
\setlength {\unitlength}{1.4cm}
\begin{picture}(5,4)
\put(1,1){\line(1,0){3}}
\put(1,1){\line(0,1){2}}
\put(1,3){\line(1,0){3}}
\put(0.7,3){$1$}
\put(0.5,1){$-1$}
\put(4.1,1){$-2$}
\put(4.1,3){$2$}
\put(1,1){\circle*{0.2}}
\put(4,1){\circle*{0.2}}
\put(1,3){\circle*{0.2}}
\put(4,3){\circle*{0.2}}
\end{picture}
\end{example}
Since the Laplace operator $\Delta_\mu$ is defined by:
\[
(\Delta_\mu f)(x) = f(x) - \displaystyle\sum_{x\sim y}f(y)\mu(x,y) ,
\]
and
\begin{equation}
 \mu (x,y) = \left\{
\begin  {matrix}
{\displaystyle\frac{1}{ deg(x)} }& \quad if (x,y) \in S\times S \\[0.7em]
0 & \quad otherwise 
\end {matrix}\right.
\end{equation}

Then the matrix of the Laplace operator $\Delta_\mu$ with respect to the basis \{ $\delta_s | 
s \in S$~\} is the following matrix:
\begin{equation}
A = \left(
\begin{matrix}
1 & -\displaystyle\frac{1}{2} & -\displaystyle\frac{1}{2} & 0\\[0.7em]
-1 & 1 & 0 & 0 \\[0.7em]
-\displaystyle\frac{1}{2} & 0 & 1 & -\displaystyle\frac{1}{2}\\[0.7em]
0 & 0 & -1 & 1
\end{matrix}\right)
\end{equation}

Then $ det( A - \alpha I_4 ) = (1 - \alpha )^2 [ (1 - \alpha)^2 - \displaystyle\frac{3}{4}] - 
\displaystyle\frac{1}{2}(1 - \alpha)^2 + \displaystyle\frac{1}{4} = 0 $\\
$\Rightarrow \alpha \in \{ 0, \displaystyle\frac{1}{2}, \displaystyle\frac{3}{2}, 2 \} 
\Rightarrow \lambda_1(\mathcal{G}(S)) =  \displaystyle\frac{1}{2}$. 

But $\mathbb{Z}$ does not have Property (T).( see \cite{bhv})

\begin{example}
The group $SL_2(\mathbb{Z})$ is generated by the matrices  $A = \left(
\begin{matrix}
1 & 1\\[0.7em]
0 & 1
\end{matrix}\right)$ and $B = \left(
\begin{matrix}
0 & -1\\[0.7em]
1 & 0
\end{matrix}\right)$.

 We consider the following generating set of the group $SL_2(\mathbb{Z})$:\\
$S = \{ -I, A, B, -A, -B, A^{-1}, B^{-1}, -A^{-1}, -B^{-1} \}$ .

The graph $\mathcal{G}(S)$ is:\\
\setlength {\unitlength}{0.8cm}
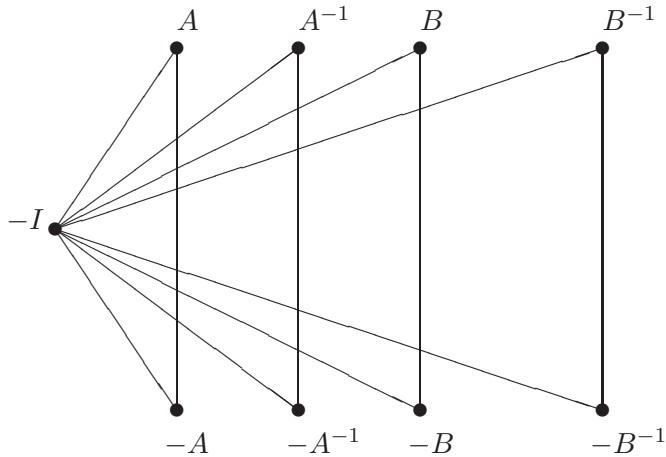
\begin{picture}(11,9)
 
\put(1,4){\line(2,3){2}}
\put(1,4){\line(4,3){4}}
\put(1,4){\line(2,1){6}}
\put(1,4){\line(3,1){9}}
\put(1,4){\line(2,-3){2}}
\put(1,4){\line(4,-3){4}}
\put(1,4){\line(2,-1){6}}
\put(1,4){\line(3,-1){9}}
\put(3,1){\line(0,1){6}}
\put(5,1){\line(0,1){6}}
\put(7,1){\line(0,1){6}}
\put(10,1){\line(0,1){6}}
\put(1,4){\circle*{0.2}}
\put(3,1){\circle*{0.2}}
\put(3,7){\circle*{0.2}}
\put(5,1){\circle*{0.2}}
\put(5,7){\circle*{0.2}}
\put(7,1){\circle*{0.2}}
\put(7,7){\circle*{0.2}}
\put(10,1){\circle*{0.2}}
\put(10,7){\circle*{0.2}}
\put(0.2,4){$-I$}
\put(2.8,0.3){$-A$}
\put(4.8,0.3){$-A^{-1}$}
\put(6.8,0.3){$-B$}
\put(9.8,0.3){$-B^{-1}$}
\put(3,7.3){$A$}
\put(5,7.3){$A^{-1}$}
\put(7,7.3){$B$}
\put(10,7.3){$B^{-1}$}
\end{picture}
\end{example}

Then the matrix of Laplace operator$\Delta_\mu$ with respect to the basis
$\{\delta_s | s \in S\}$ is the following matrix:\\
\begin{equation}
A = \left(
\begin{matrix}
1 & \displaystyle\frac{1}{8} & \displaystyle\frac{1}{8} & \displaystyle\frac{1}{8} &
\displaystyle \frac{1}{8} & \displaystyle\frac{1}{8 }&\displaystyle\frac{1}{8}
&\displaystyle\frac{1}{8} & \displaystyle\frac{1}{8}\\[0.7em]
\displaystyle\frac{1}{2} & 1 &\displaystyle\frac{1}{2} & 0 & 0 & 0 & 0 & 0 & 0\\[0.7em]
\displaystyle\frac{1}{2}  &\displaystyle\frac{1}{2} & 1& 0 & 0 & 0 & 0 & 0 & 0\\[0.7em]
\displaystyle\frac{1}{2} & 0 & 0 & 1 &\displaystyle\frac{1}{2} & 0 & 0 & 0 & 0 \\[0.7em]
\displaystyle\frac{1}{2} & 0 & 0 &\displaystyle\frac{1}{2}& 1 & 0 & 0 & 0 & 0 \\[0.7em]
\displaystyle\frac{1}{2} & 0 & 0 & 0 & 0 & 1 &\displaystyle\frac{1}{2} & 0 & 0  \\[0.7em]
\displaystyle\frac{1}{2} & 0 & 0 & 0 & 0  &\displaystyle\frac{1}{2} & 1 & 0 & 0  \\[0.7em]
\displaystyle\frac{1}{2} & 0 & 0 & 0 & 0 & 0 & 0 & 1 &\displaystyle\frac{1}{2} \\[0.7em]
\displaystyle\frac{1}{2} & 0 & 0 & 0 & 0 & 0 & 0  &\displaystyle\frac{1}{2} & 1 
\end{matrix}\right)
\end{equation} 

Computing det($A-\alpha I_9$) = $[ (1-\alpha )^2-\displaystyle\frac{1}{4}]^3(\displaystyle\frac{3}{2}-\alpha)
(\alpha^2-\displaystyle\frac{5}{2}\alpha) = 0 \Rightarrow $\\
$\Rightarrow \alpha \in \{ 0, \displaystyle\frac{1}{2}, \displaystyle\frac{3}{2},
\displaystyle\frac{5}{2} \} \Rightarrow \lambda_1( \mathcal{G}(S)) = \displaystyle\frac{1}{2}$.

But $SL_2(\mathbb{Z})$ does not have Property (T). (see \cite{bhv})

\textbf{These two examples shows that $\displaystyle\frac{1}{2}$ is the best constant in Zuk's criterion and cannot be improved.}

{\footnotesize University of Bucharest, Romania }

{\footnotesize e-mail: traianpr@yahoo.com}

\end{document}